\documentclass[a4,reqno]{amsart}
\usepackage{amsmath,amssymb,amsthm,enumerate,bm}
\usepackage[monochrome]{color}

\newcommand{\Q}{{\mathbb Q}}
\newcommand{\Z}{{\mathbb Z}}
\newcommand{\C}{{\mathbb C}}
\newcommand{\N}{{\mathbb N}}
\newcommand{\Ex}{\lambda}
\renewcommand{\i}{{\bm i}}

\theoremstyle{plain}
\newtheorem{thm}{Theorem}
\newtheorem*{thm*}{Theorem}

\newtheorem*{defi*}{Definition}
\newtheorem{rem}{Remark}

\newtheorem{cor}{Corollary}
\newtheorem{conj}{Conjecture}
\newtheorem{prop}{Proposition}
\newtheorem{ex}{Example}

\newcommand{\wi}{\widetilde}

\newcommand{\bt}{\mathbf{t}}
\newcommand{\Ph}{q}

\begin{document}
\title{Curious congruences for cyclotomic polynomials}
\keywords{Cyclotomic polynomials, 
Euler's totient function, Jordan totient function, 
congruence relation}
\subjclass[2020]{Primary~11A25,~Secondary~11A07,~11R18}
\maketitle
\centerline{Shigeki Akiyama and Hajime Kaneko}
\hspace*{1cm}

\begin{center}
{\footnotesize
Institute of Mathematics / Research Core for Mathematical Sciences, 
University of Tsukuba, Tsukuba, JAPAN. 
(e-mail: akiyama@math.tsukuba.ac.jp, kanekoha@math.tsukuba.ac.jp)\\
Data openly available in a public repository: https://arxiv.org/abs/2204.11267
}
\end{center}
\begin{abstract}
Let $\Phi_n^{(k)}(x)$ be the $k$-the derivative of the $n$-th cyclotomic
polynomial. We are interested in the values $\Phi_n^{(k)}(1)$ for fixed 
positive integers $n$. 
\par
D.~H.~Lehmer proved that $\Phi_n^{(k)}(1)/
\Phi_n(1)$ is a polynomial of the Euler totient function $\phi(n)$ and
the Jordan totient functions and gave its explicit formula.
In this paper, we give a quick proof that
$\Phi_n^{(k)}(1)/\Phi_n(1)$ is a polynomial of them 
without giving the explicit form.
\par
In the final section, we deduce some curious congruences: $2\Phi^{(3)}_n(1)$
is divisible by $\phi(n)-2$. Moreover, if $k$ is greater than 1, then
$\Phi^{(2k+1)}_n(1)$ is divisible by $\phi(n)-2k$.
The proof depends on a new combinatorial identity for general self-reciprocal polynomials over $\Z$, 
which gives rise to
a formula that expresses the value
$\Phi_n^{(k)}(1)$ as a
$\Z$-linear combination of the coefficients in the minimal polynomial
of $2\cos(2\pi/n)-2$.

As a supplement, we show the monotonic increasing property of $\Phi_n(x)$ on $[1,\infty)$ in two ways.
\end{abstract}

\section{Introduction}\label{section1}

The $n$-th cyclotomic polynomial

\begin{equation}
\label{Def0}
\Phi_n(x)=\prod_{\stackrel{0<d<n}{(d,n)=1}} \left(x-\exp\left(\frac{2\pi d \i}{n}\right)\right)
\end{equation}
is the minimal polynomial of the $n$-th primitive roots of unity over $\Q$. 
It is an irreducible polynomial in $\Z[x]$ of degree $\phi(n)$
where $\phi$ is the Euler totient function. 
From the relation $x^n-1=\prod_{d\mid n} \Phi_d(x)$, the well known formula

\begin{equation}
\label{Def1}
\Phi_n(x)=\prod_{d \mid n} (x^d-1)^{\mu(n/d)}
\end{equation}
is derived by M\"obius inversion. Here $\mu$ is the M\"obius function. 
Motose \cite{Motose:04, Motose:05} surmised that $\Phi_n(x)$ is an increasing function\footnote{He did not give a proof of this fact,
see the sentence before Theorem 3 in \cite{Motose:04}.} for $x>1$.
We start with a simple proof of this fact.
It is probably known but we did not find it in the
literature.

\begin{thm}
\label{Motose}
\textcolor{blue}{
For $j=1,\dots, \phi(n)$ we have}
\label{Pos}
$$
\Phi_{n}^{(j)}(1)>0.
$$
\textcolor{blue}{Consequently} $\Phi_n^{(k)}(x)$ is strictly
increasing for $x\ge 1$
and $k=0,1,\dots, \phi(n)-1.$
\end{thm}

\proof
Since $\Phi_1^{(1)}(1)=1$ and
$\Phi_2^{(1)}(1)=1$, 
we may assume that $n\ge 3$. Then we have
\begin{align*}
\Phi_n(x)
&=\prod_{\stackrel{0<d<n/2}{(d,n)=1}} 
\left(x-\exp\left(\frac{2\pi d \i}{n}\right)\right)
\left(x-\exp\left(\frac{2\pi(n-d) \i}{n}\right)\right)\\
&=\prod_{\stackrel{0<d<n/2}{(d,n)=1}} 
\left(x^2 -2 \cos \left(\frac{2\pi d}n\right) x +1\right).
\end{align*}
Since all coefficients of 
$$
(x+1)^2+b(x+1)+1=x^2+ (b+2) x + b+2
$$
with $b\in (-2,2)$ are positive, the expansion $\Phi_n(x+1)=
\sum_{j=0}^{d} \Phi_n^{(j)}(1)/j! \cdot x^j$ at $x=0$ 
have positive coefficients $\Phi_n^{(j)}(1)/j!$ for $j\le d=\phi(n)$. 
This proves the theorem.
\qed
\bigskip

\begin{rem}
\label{Opt}
The inequality $x\ge 1$ in Theorem \ref{Pos} is sharp.
If $p$ is an odd prime, then
$$
\Phi_{2p}(x)=\frac{1-(-x)^p}{1+x}.
$$
It is easy to confirm
\begin{equation}
\label{Neg}
\Phi_{2p}'\left(1-\frac{1}{\sqrt{p}}\right)
=\frac{\left(2p-\sqrt{p}+\frac{1}{\sqrt{p}}-1\right)\left(1-\frac{1}{\sqrt{p}}\right)^{p-1}-1}{\left(2-\frac{1}{\sqrt{p}}\right)^2}<0.
\end{equation}
Thus there exists no $\varepsilon>0$ that 
$\Phi_{n}(x)$ is increasing on $[1-\varepsilon,\infty)$ for any $n\geq 1$. 
\end{rem}

\begin{rem}
\label{Alt}
There is an alternative proof that only works for $k=0$, 
giving a starting point for this paper.
Since $\Phi_1(x)=x-1$, we may assume that $n\ge 2$. From 
\begin{equation}
\label{MZ}
\sum_{d\mid n} \mu(d)=0
\end{equation}
for $n\ge 2$, we may replace $x^d-1$ by $(x^d-1)/(x-1)$ in (\ref{Def1}). 
Taking $x\rightarrow 1$, we have
$$
\Phi_n(1)=\prod_{d\mid n} d^{\mu(n/d)}
$$
which is a positive integer. We see this is rewritten as
$\Phi_n(1)=\exp(\Lambda(n))$ with the von Mangoldt function
$$
\Lambda(n):=\begin{cases} \log p & n=p^e\ (p\text{ prime})\\ 
                                           0   & \text{otherwise},
\end{cases}
$$
which plays a crucial role in analytic number theory. 
The fact above was proved by Lebesgue \cite{Lebesque:859}. 
Using (\ref{MZ}), we also have
\begin{align}
\log \Phi_n(x)&= \sum_{d\mid n} \mu\left(\frac nd\right) \log\left(\frac{x^d-1}{x-1}\right)\notag\\
\frac{\Phi_n'(x)}{\Phi_n(x)}
&=\sum_{d\mid n} \mu\left(\frac nd\right) 
\left(\frac{(d-1)x^{d-2}+(d-2)x^{d-3}+\dots +1}{x^{d-1}+x^{d-2}\dots +1}\right).\label{D}
\end{align}
Letting $x\to 1$ and using (\ref{MZ}) again, we obtain
$$
\frac{\Phi_n'(1)}{\Phi_n(1)}
=\sum_{d\mid n} \mu\left(\frac nd\right) 
\frac{d-1}2=\frac 12\sum_{d\mid n} \mu\left(\frac nd\right)d=\frac{\phi(n)}2.
$$
Thus we see that 
\begin{equation}
\label{Diff}
\Phi_n'(1)=\frac 12 \phi(n)\Phi_n(1)\ge 1>0,
\end{equation}
which was proved by H\"{o}lder \cite{Holder:36} (c.f. \cite[Lemma 10]{BHM}).
Now we consider $\Phi_n(z)$ as a polynomial of complex variable $z\in \C$. 
Recalling Gauss--Lucas theorem, any root of $\Phi_n'(z)$ lies 
in the convex hull of the roots of $\Phi_n(z)$ in the complex plane. 
Therefore from (\ref{Def0})
and $n\ge 2$, the 
real function $\Phi_n'(x)$ has no root in $x\ge 1$. 
This implies $\Phi_n'(x)>0$ for $x\ge 1$ since $\Phi_n'$ is continuous.
\end{rem}

\begin{rem}
Let $p$ be an odd prime. Then
(\ref{Neg}) and 
(\ref{Diff}) imply that there exists a real root of 
$\Phi'_{2p}(x)$ in the interval $(1-1/\sqrt{p},1)$.
\end{rem}

Jordan totient function is defined by $J_k(n)=\sum_{d\mid n} \mu(n/d)d^k$.
This is multiplicative and we have $$
J_k(n)=n^k\prod_{p\mid n} \left(1-\frac 1{p^k}\right)
$$
where $p$ runs over prime divisors of $n$. 
Clearly  $J_k(n)$ is a generalization of 
the Euler totient function $\phi(n)=J_1(n)$.  
The name came from C.~Jordan who studied linear groups over $\Z/n\Z$
and deduced, e.g., 
$$
\mathrm{Card}(GL_k(\Z/n\Z))= n^{\frac{k(k-1)}2} \prod_{j=1}^k J_j(n).
$$

As we observed in Remark \ref{Alt}, the special values
$\Phi_n^{(k)}(1)$ give important arithmetic functions such as the
von Mangoldt function and the Euler totient function. 
Lehmer \cite{Lehmer:66} gave an explicit formula of  
$\Phi_n^{(k)}(1)/\Phi_n(1)$ 
as a polynomial of $\phi(n)$ and $J_{2j}(n)$ over $\Q$, using Stirling numbers and Bernoulli numbers, 
see \cite{SC, HM:21, MESS:20, Sanna} for related developments. 
Here we give a quick proof of this fact but without the explicit form of the polynomial.

\begin{thm} [\cite{Lehmer:66}]
For $n\ge 2$, 
\label{JordanExpansion}
$\Phi_n^{(\ell)}(1)/\Phi_n(1)$ is expressed as a polynomial
of $\phi(n)$ and $J_{2j}(n)\ (1\le j \le (\ell+1)/2)$ over $\Q$, 
and
its value is a positive integer\footnote{Clearly $\Phi_n^{(\ell)}(x)=0$
for $\phi(n)<\ell$.} for $\phi(n)\ge \ell$.
\end{thm}

\begin{proof}
Applying Leibniz formula to (\ref{D}),
$$
\frac{\Phi_n^{(k+1)}(x)}{\Phi_n(x)}
= \sum_{\ell=0}^{k} \binom{k}{\ell} \frac{\Phi_n^{(\ell)}(x)}{\Phi_n(x)} \sum_{d\mid n}
\mu\left(\frac nd\right)
\frac{\partial^{k-\ell}}{\partial x^{k-\ell}}
\left(\frac{1-d x^{d-1}+d x^d-x^d}{(x-1) (x^d-1)}\right).
$$
Substituting $x$ by $1+t$, we get the Taylor expansion at $t=0$:
\begin{equation}
\label{Taylor}
\frac{((d-1) t-1) (t+1)^{d-1}+1}{t \left((t+1)^d-1\right)}
=
\frac{d-1}{2}+\frac{d^2-6 d+5}{12} t+\frac{-d^2+4 d-3}{8}
   t^2+O\left(t^3\right).
\end{equation}
Regarding $d$ as a real variable,
we see that the numerator of (\ref{Taylor}) has the expansion at $t=0$
of the form
$$
\sum_{j\ge 0} d f_j t^{j+2} \quad f_j\in \Q[d],\ \mathrm{deg}(f_j)=j+1,\ f_0=(d-1)/2.
$$
Similarly, the denominator has the form
$$
\sum_{j\ge 0} d g_j t^{j+2} 
\quad g_j\in \Q[d],\ \mathrm{deg}(g_j)=j,\ g_0=1.
$$
Thus, the $\ell$-th Taylor coefficient of their quotient
is a polynomial of $d$ whose degree does not exceed $\ell+1$.
Using these Taylor coefficients, we recursively obtain
the explicit formula for $\Phi_n^{(\ell)}(1)/\Phi_n(1)$.
Thus 
$\Phi_n^{(\ell)}(1)/\Phi_n(1)$
a polynomial on $J_{1}(n), J_{2}(n), \dots, J_{\ell+1}(n)$ over $\Q$.
Moreover since 
$$
\frac{((d-1) t-1) (t+1)^{d-1}+1}{t \left((t+1)^d-1\right)}
-\frac{d}{2 (t+1)}=
\frac{d}{2 (t+1)} \frac{(t+1)^d+1}{(t+1)^d-1}
-\frac{1}{t}
$$
is an even function on $d$, the terms $d^{2k+1}$ with $k=1,2,\dots$
do not show, i.e., $J_{2k+1}(n)\ (k=1,2,\dots)$ never appear.
By Theorem \ref{Pos}, $\Phi_n^{(\ell)}(1)/\Phi_n(1)>0$ for $\phi(n)\ge \ell$.
Since $$
\Phi_n(1)=\exp(\Lambda(n))=\begin{cases} p & n=p^e\ (p: \text{prime})\\
                                         1 & \text{otherwise},
\end{cases}
$$
it suffices to show 
$\Phi_{p^e}^{(\ell)}(1)\equiv 0 \pmod{p}$. By
$$
\Phi_{p^e}(x)=\frac{x^{p^e}-1}{x^{p^{e-1}}-1}=\Phi_{p}(x^{p^{e-1}}),
$$
the case $e>1$ is plain and the case $e=1$ remains to be settled.
Indeed we have,
$$\Phi_{p}^{(\ell)}(1)=\sum_{j=0}^{p-1} j(j-1)\cdots (j-\ell+1)=
\frac{p(p-1)\cdots (p-\ell)}{\ell+1} \equiv 0 \pmod{p}.$$
\end{proof}


\begin{cor}
We have 
\label{23}
\begin{align*}
\frac{\Phi_n^{(2)}(1)}{\Phi_n(1)}=&
\frac{J_2(n)}{12}+\frac{\phi(n)^2}4-\frac{\phi(n)}2,\\
\frac{\Phi_n^{(3)}(1)}{\Phi_n(1)}=&\frac{(\phi(n)-2)(J_2(n)+\phi(n)(\phi(n)-4))}8,\\
\frac{\Phi_n^{(4)}(1)}{\Phi_n(1)}=&
\frac{1}{240} \Bigg(30 J_2(n) \phi (n)^2-180 J_2(n) \phi (n)+5 J_2(n)^2+220 J_2(n)-2J_4(n)\\
&+15 \phi (n)^4-180 \phi (n)^3+660 \phi (n)^2-720 \phi (n)\Bigg), \\
\frac{\Phi_n^{(5)}(1)}{\Phi_n(1)(\phi(n)-4)}
=&\frac{1}{96} 
\Big(3 \phi(n)^4-48 \phi(n)^3+10 J_2(n) \phi(n)^2+228
   \phi(n)^2\label{5} \\
&-80 J_2(n) \phi(n)-288 \phi(n)+5 J_2(n)^2+100
   J_2(n)-2 J_4(n)\Big).
\end{align*}
\end{cor}

Let 
\begin{equation*}
\label{bn}
\Phi_n(x+1)=:\sum_{h=0}^{\phi(n)}b_n(h) x^h \ \text{ with }\ b_n(h)=\frac{1}{h!}
\Phi_n^{(h)}(1)\in \Z.
\end{equation*}
Lehmer \cite{Lehmer:66} further stated an interesting observation on
the coefficients $b_n(h)$.
For a real $R$, set $R^{[\ell]}:=R(R-1)\cdots (R-\ell+1)$. 
For a positive integer $r$, let $t_r:=J_r(n)/(2r)$. 
We define Bernoulli numbers $B_m$ ($m\geq 0$) by 
\[
\frac{t e^t}{e^{t}-1}=\sum_{n=0}^{\infty} B_n \frac{t^n}{n!}. 
\]
Under the setting above, he claimed that 
\[
\frac{b_n(h)}{\Phi_n(1)}=t_1^{[h]}+2\sum_{\ell=1}^{\infty} B_{2\ell}
\binom{h}{2\ell} (t_1-\ell)^{[h-2\ell]}\Omega_{\ell},
\]
but the general form of $\Omega_{\ell}$ is not given. He only wrote the first few terms:
\begin{align*}
\Omega_1&=t_2, \\
\Omega_2&=t_4-5t_2^{[2]}, \\
\Omega_3&=t_6-7t_4(t_2-1)+\frac{35}{3} t_2^{[3]}+\frac{14}{3}t_2,\\
\Omega_4&=t_8-\frac{20}{3}t_6(t_2-1)-\frac{7}{3}t_4^{[2]}
+\frac{70}{3} t_4(t_2-1)^{[2]}\\
&\hspace{20mm}-\frac{175}{9}t_2^{[4]}
+\frac{10}{3}t_6-\frac{280}{9}t_2^{[2]}+\frac{290}{9}t_2.
\end{align*}

Both Corollary \ref{23} and this observation suggest the following: 

\begin{conj}
\label{Div}
For any non-negative integer $k$, 
$\Phi_n^{(2k+1)}(1)/\Phi_n(1)$ is divisible by $\phi(n)-2k$ in the polynomial
ring 
$\Q\left(\phi(n), J_{2}(n), J_4(n),\dots, J_{2\lfloor \ell+1)/2\rfloor}(n)\right)$. 
\end{conj}

We checked its validity for $k\le 15$.
The goal of this paper is to prove intimately related divisibility:
$$
\Phi_n^{(2k+1)}(1)
\ \text{is divisible by}\ \phi(n)-2k\ \text{ in } \Z ,
$$
for $k\ge 1$, see Theorem \ref{cor:2-1}. 
(The dividend should be doubled for the case $k=1$.) 
We did not find yet a special meaning for this divisibility.
For a fixed $n$, such divisibility is proved 
using Theorem \ref{JordanExpansion}
\textcolor{blue}{and Proposition \ref{Triv} below.
However, such an individual proof does not seem to extend to the general case.}
Note that $\Phi_n^{(2k+1)}(1)/\Phi_n(1)$ is likely to be divisible by $\phi(n)-2k$ for $k>1$
but it can not be shown by our method.

Let $\Ex(m)$ be the Carmichael lambda function, i.e., 
the exponent of $(\Z/m\Z)^*$, the unit group of the ring
$\Z/m\Z$ (c.f. \cite{Car,SC}). For an odd prime $p$,
$\Ex(p^e)=\phi(p^e)$ holds since $(\Z/p^e\Z)^*$ is cyclic.
From
$$
(\Z/2^e\Z)^*\simeq \Z/2\Z \oplus \Z/2^{e-2}\Z
$$
for $e\ge 2$, we have
$$
\Ex(2^e)= \begin{cases} 
1 & e=1\\
2 & e=2\\
2^{e-2} & e\ge 3.
\end{cases}
$$
For a prime $p$ and a positive integer $e$, we write 
$\Ex(p^e)\parallel k$ if both $\Ex(p^e)\mid k$ and $\Ex(p^{e+1})\nmid k$ hold.
Proposition \ref{Triv} may be
known, but we give a proof for self-containedness.

\begin{prop}[Trivial congruence]
\label{Triv}
For $k\ge 3$ and $n\ge k+2$, we have
$$
J_{k}(n)\equiv 0 \pmod{\prod_{\Ex(p^e)\parallel k} p^e}.
$$
For $M>\prod_{\Ex(p^e)\parallel k} p^e$
and any $n_0\in \N$, there exists $n\ge n_0$ such that 
$
J_{k}(n)\not \equiv 0 \pmod{M}.
$
\end{prop}

\begin{proof}
There are only finitely many prime $p$ such that
$\Ex(p^e) \parallel k$.
For a prime factor $q$ of $n$, $q^k-1$ is a factor of $J_k(n)$.
The condition $\Ex(p^e) \mid k$ implies 
$q^k-1\equiv 0 \pmod{p^e}$ for each $q$ which is coprime
with $p$. 
Assume that 
\begin{equation}
\label{Suff}
n>\max \{ p \ :\ \Ex(p^e)\mid k \}.
\end{equation}
If $n$ has two distinct prime factors $p_1$ and $p_2$, then 
$J_{k}(n)$ is divisible by $(p_1^k-1)(p_2^k-1)$. We see
$p_j^k-1$ is divisible by $p^e$ with
$p\neq p_j$
if $\Ex(p^e)\mid k$. This implies that 
$(p_1^k-1)(p_2^k-1)$ is divisible by
$\prod_{\Ex(p^e)\parallel k} p^e$. 
Thus we may assume that $n$ is a power of a prime $q$ and 
$\Ex(q^e)\mid k$, i.e., $n=q^{\ell}$ and $\ell\ge 2$.
In this case, $J_k(n)$ is divisible by $q^{\ell k}-q^{(\ell-1)k}=
q^{(\ell-1)k}(q^k-1)$. We see
$$
\prod_{\Ex(p^e)\parallel k\ p\neq q} p^e\ \mid \ q^k-1.
$$
When $q\neq 2$, since $e\le \Ex(q^e)\le k$, we see 
$q^{e}$ divides $q^{(\ell-1)k}$ and 
the required congruence holds. 
For $q=2$ we only have $e-1\le \Ex(2^e)\le k$ and hence $e\le 2k$.
So additionally if $\ell>2$, then $q^e \mid q^{(\ell-1)k}$ holds.
Therefore our discussion fails only when  $n=2^2$, $\Ex(2^e)\mid k$ and $e>k$. 
This happens when $2^{k-1}\le k$, that is, $k\le 2$.
Summing up if $k\ge 3$,  (\ref{Suff}) implies our congruence. 
Moreover (\ref{Suff}) holds if $n>k+1$, because the worst case happens
when $k+1$ is an odd prime.

Take $M>\prod_{\Ex(p^e)\parallel k} p^e$ and any $n_0\in \N$. 
There exists a prime power factor $p^{e+1}$ of $M$ 
that $\Ex(p^{e+1})$ does not divide $k$. 
From the definition of the exponent, 
there exists $t\in \N$ which is coprime to $p$ that
$t^k\not \equiv 1 \mod{p^{e+1}}$. By Dirichlet's theorem,
there exists a prime $q\ge n_0$ that $q\equiv t \pmod{p^{e+1}}$. Then $J_{k}(q)=q^k-1
\not \equiv 0 \mod{M}$.
\end{proof}

Here is a table of the first numbers appearing in Proposition \ref{Triv}, see also \cite{OEIS}.

\begin{table}[htbp]
  \centering
  \begin{tabular}{c|c|c|c|c|c|c|c|c|c|c|c|}
$k$& odd & 2 & 4 & 6 & 8 & 10& 12 & 14 & 16 & 18 & 20 
\\ \hline
$
\prod_{\Ex(p^e)\parallel k} p^e
$ & 2 &24 & 240& 504 & 480 & 264 & 65520&  24& 16320& 28728& 13200
  \end{tabular}
\end{table}

\section{Congruences for self-reciprocal polynomials}\label{section2}

Let $\Ph$ be a positive integer and 
$t_j$ be a complex number for $1\le j\le \Ph$.
Let
\begin{align}\label{eqn:2-1}
\sum_{h=0}^{2\Ph}b(\bt;h) x^h:=
\prod_{j=1}^{\Ph}(x^2+t_jx+t_j)
\end{align}
and 
\begin{align}\label{eqn:2-2}
\sum_{\ell=0}^{\Ph} a(\bt;\ell) y^{\ell}:=
\prod_{j=1}^{\Ph}(y+t_j)
\end{align}
where $\bt=(t_j)_{1\leq j\leq \Ph}$.
Our key result is a special 
combinatorial equality between $b(\bt;h)$ and $a(\bt;\ell)$.

\begin{prop}\label{thm:main}
Let $\bt=(t_j)_{1\leq j\leq \Ph}$ be a sequence of complex numbers. 
Then, for any $h$ with $0\leq h\leq 2\Ph$, we have 
\begin{align*}
b(\bt;h)=\sum_{\ell=\max\{0,h-q\}}^{\lfloor h/2\rfloor}
\binom{\Ph-\ell}{h-2\ell} a(\bt;\ell). 
\end{align*}
\end{prop}

\begin{proof} We have 
\begin{align*}
\prod_{j=1}^q (x^2+t_j x + t_j) &= (x+1)^q \prod_{j=1}^q \left(\frac {x^2}{x+1}+ t_j\right)\\
&= (x+1)^q \sum_{\ell=0}^q a(\bt;\ell) \left(\frac {x^2}{x+1}\right)^\ell\\
&= \sum_{\ell=0}^q a(\bt;\ell) x^{2\ell} (x+1)^{q-\ell}\\
&= \sum_{\ell=0}^q a(\bt;\ell) x^{2\ell} \sum_{h=0}^{q-\ell} \binom{q-\ell}{h} x^{h} \\
&= \sum_{\ell=0}^q a(\bt;\ell) \sum_{h=2\ell}^{q+\ell} \binom{q-\ell}{h-2\ell} x^{h}\\
&= \sum_{h=0}^{2q} \left(\sum_{\ell=\max\{0, h-q\}}^{\lfloor h/2\rfloor} 
\binom{q-\ell}{h-2\ell} a(\bt;\ell)\right) x^{h}.
\end{align*}
Comparing the coefficients, we get the desired result. 
\end{proof}


\begin{rem}\label{rem:index}
Under the usual convention $\binom{m}{n}=0$ for non negative integers $m,n$,
 with $m<n$,
Proposition \ref{thm:main} is rephrazed as
\begin{align*}
b(\bt;h)=\sum_{\ell=0}^{\lfloor h/2\rfloor}
\binom{\Ph-\ell}{h-2\ell} a(\bt;\ell). 
\end{align*}
\end{rem}

\begin{ex}
Let $a$ be a complex number. 
Setting $t_j=a$ for $j=1,\dots,\Ph$, we obtain
$$(x^2+ax+a)^{\Ph}=
\sum_{h=0}^{2\Ph}  \left(\sum_{\ell=0}^{\lfloor h/2\rfloor}
\binom{q-\ell}{h-2\ell} \binom{\Ph}{\ell}a^{\Ph-\ell}\right) x^h.$$
This is also shown directly by the binomial theorem. 
\end{ex}

\begin{ex}
Setting $t_j=\pm j$ for $j=1,\dots,\Ph$, 
we have
\begin{align*}
\prod_{j=1}^{\Ph}(x^2-jx-j)&=\sum_{h=0}^{2\Ph}
\left(\sum_{\ell=0}^{\lfloor h/2\rfloor}
\binom{\Ph-\ell}{h-2\ell}  s(q+1,\ell+1) \right) x^h,\\
\prod_{j=1}^{\Ph}(x^2+jx+j)&=\sum_{h=0}^{2\Ph}
\left(\sum_{\ell=0}^{\lfloor h/2\rfloor}
\binom{\Ph-\ell}{h-2\ell}  |s(q+1,\ell+1)| \right) x^h
\end{align*}
where $s(m,n)$ is the Stirling number of the 1-st kind
defined by
$$
x(x-1)\cdots (x-m+1)=\sum_{n=0}^m s(m,n) x^n.
$$
\end{ex}


We apply Proposition \ref{thm:main} to general self-reciprocal polynomials of even degree.
Let $f(x)\in \Z[x]$ be a polynomial with degree $2q$ ($q\geq 1$). 
Suppose that $f(x)$ is self-reciprocal, that is, $x^{2q}f(x^{-1})=f(x)$. 
It is easily seen by induction on $q$ that there exists $g(y)\in \Z[y]$ such that 
$f(x)=x^q g(y)$ with $y=x+x^{-1}$. Let 
$$
f(x+1):=\sum_{h=0}^{2q} \beta(h) x^h, \quad
g(y+2):=\sum_{\ell=0}^q \alpha(\ell)y^{\ell},
$$
where $\beta(2q)=\alpha(q)$. 
Then we have the following: 
\begin{prop}\label{recip_1}
For any $h$ with $0\leq h\leq 2q$, we have 
$$
\beta(h)=\sum_{\ell=0}^{\lfloor h/2\rfloor} \binom{\Ph-\ell}{h-2\ell} \alpha(\ell).
$$
\end{prop}
\begin{proof}
We now reduce Proposition \ref{recip_1} to 
Proposition \ref{thm:main}. 
Let 
\[
g(y)=:\beta(2q)\prod_{j=1}^q (y+\gamma_j),
\]
where $\gamma_1,\ldots,\gamma_q$ are complex numbers. 
Then we see 
\[
f(x)=\beta(2q)\prod_{j=1}^q (x^2+\gamma_j x+1),
\]
and so 
\[
f(x+1)=\beta(2q)\prod_{j=1}^q \big(x^2+(\gamma_j+2) x+\gamma_j+2\big).
\]
On the other hand, using 
\[
f(x)=x^q\cdot \beta(2q)\prod_{j=1}^q\left(x+\frac1x+\gamma_j\right)
=x^q g(y),
\]
we see 
\[
g(y+2)=\beta(2q)\prod_{j=1}^q\left(y+\gamma_j+2\right).
\]
Hence, (\ref{eqn:2-1}) and (\ref{eqn:2-2}) hold with $t_j=\gamma_j+2$ for $j=1,2,\ldots, q$. 
Therefore, Proposition \ref{recip_1} follows from 
Proposition \ref{thm:main} and Remark \ref{rem:index}. 
\end{proof}
Proposition \ref{recip_1} leads to the following congruences.
\begin{prop}\label{recip_2}
\begin{enumerate}[(i)]
\item $2f'''(1)$ is divisible by $2q-2$. Moreover, if $q$ is even, then $f'''(1)$ is divisible by $2q-2$. 
\item Suppose that $k\geq 2$. Then $f^{(2k+1)}(1)$ is divisible by 
$2q-2k$. 
\end{enumerate}
\end{prop}

\begin{proof}
For the proof of (i), we may assume that 
$2q\geq 3$. Since $2f'''(1)=12 \beta(3)$, (i)
follows from Proposition \ref{recip_1} and 
\begin{align*}
12 \binom{q}{3}&=(2q-2)\cdot 
q(q-2). 
\end{align*}
The latter part of (i) is similarly proved 
because $q$ is even. \par
For the proof of (ii), we may assume that 
$2q\geq 2k+1$. Then (ii) follows from 
$f^{(2k+1)}(1)=(2k+1)!\beta(2k+1)$ and 
Proposition \ref{recip_1}. 
Set 
\[
\wi{c(\ell)}:=(2k+1)! \binom{q-\ell}{2k+1-2\ell}.
\]
We shall show that $\wi{c(\ell)}$ is divisible by $2q-2k$. We may assume that $\wi{c(\ell)}\ne 0$. 
Note that $(q-2k+1)(q-2k)$ 
is even. 
If $\ell=0$, then we see by $k\geq 2$ that 
\begin{align*}
\wi{c(0)}=q(q-1)\cdots (q-k) \cdots (q-2k+1)(q-2k)
\end{align*}
is divisible by $2q-2k$. 
Moreover, if $1\leq \ell\leq k$, then 
\[
\wi{c(\ell)}=
\frac{(2k+1)!}{(2k+1-2\ell)!} \cdot
(q-\ell) \cdots (q-k) \cdots (q-2k+\ell)
\]
is divisible by $2q-2k$ because $(2k+1)!/(2k+1-2\ell)!$ is even. 
\end{proof}

Recall that $\Phi_n(x)$ is self-reciprocal.
For $f(x)=\Phi_n(x)$, the corresponding $g(y)$ is the minimum polynomial of $\zeta_n+\zeta_n^{-1}$, where $\zeta_n$ is a primitive $n$-th root of unity. Proposition \ref{recip_1} gives a formula to express
$\Phi_n^{(k)}(1)$ as a $\Z$-linear combination of the coefficients of
the minimal polynomial of $\zeta_n+\zeta_n^{-1}-2$.
 Proposition \ref{recip_2} includes our curious congruences on cyclotomic polynomials.

\begin{thm}\label{cor:2-1}
\begin{enumerate}[(i)]
\item $2\Phi_n'''(1)$ is divisible by $\phi(n)-2$. In particular, if $\phi(n)$ is divisible by $4$, then $\Phi_n'''(1)$ is divisible by $\phi(n)-2$. 
\item Suppose that $k\geq 2$. Then $\Phi_n^{(2k+1)}(1)$ is divisible by 
$\phi(n)-2k$. 
\end{enumerate}
\end{thm}

\section*{Acknowledgments}

We would like to thank Pieter Moree and Michel Marcus
for their comments and references to the earlier version. The presentation of 
this paper is largely improved by 
the suggestions of anonymous referees. In particular, 
we could reach the current concise form of Proposition \ref{thm:main}
by a suggestion asking for a combinatorial reformulation.
This research was partially supported by JSPS grants (20K03528, 17K05159, 21H00989, 19K03439). 
The authors declare no conflicts of interest associated with this manuscript. 

{\bf Note added in proof.} After our presentation at RIMS on 12 Oct 2022,
T.~Matsusaka informed us of a proof of Conjecture \ref{Div}. G.~Shibukawa
told us that this proof of Theorem \ref{Motose} is published later
in Japanese textbooks for cyclotomic polynomials by K.~Motose. We also got to know that G.~Shibukawa introduced a similar method to Proposition \ref{thm:main} in Fibonacci Quart. 58 (2020), no. 5, 200-221. We are hoping to discuss these in a future work.

\end{document}